\documentclass[12pt,a4paper]{article}
\usepackage[english]{}
\usepackage{amsfonts}
\usepackage{amsthm}
\usepackage[ansinew]{inputenc}
\usepackage{eufrak}
\usepackage[sumlimits,intlimits, namelimits]{amsmath}
\usepackage{graphicx}
\usepackage{flafter}
\usepackage{color}
\usepackage{multicol}
\usepackage{fancybox}
\usepackage{rotating}
\usepackage{amssymb}

\def\N{{\mathbb N}}
\def\Z{{\mathbb Z}}
\def\R{{\mathbb R}}

\def\p{\partial}
\def\bs{\backslash}
\def\weakstar{\stackrel{\ast}{\rightharpoonup}}

\def\mcL{\mathcal{L}}

\def\calL{\mathcal{L}}
\def\calM{\mathcal{M}}
\def\calS{\mathcal{S}}
\def\calN{\mathcal{N}}
\def\calV{\mathcal{V}}
\def\calA{\mathcal{P}}
\def\calH{\mathcal{H}}

\def\diam{\mbox{diam}\,}

\def\Eloc{\mbox{$E_{\ell oc}$}}

\def\mutildeN{\tilde{\mu}_N}
\def\lambdatildeN{\tilde{\lambda}_N}
\def\nutildeN{\tilde{\nu}_N}
\def\muttildeN{\tilde{\tilde{\mu}}_N}
\def\lambdattildeN{\tilde{\tilde{\lambda}}_N}
\def\nuttildeN{\tilde{\tilde{\nu}}_N}

\newcommand{\be}{\begin{equation}}
\newcommand{\ee}{\end{equation}}
\newcommand{\Proof}{\par\noindent{\bf Proof. }}
\newcommand{\eop}{\nopagebreak\hspace*{\fill}$\Box$\smallskip}

\newtheorem{theorem}{Theorem}[section]
\newtheorem{proposition}[theorem]{Proposition}

\newtheorem{lemma}[theorem]{Lemma}

\theoremstyle{definition}

\definecolor{grau}{gray}{0.80}

\begin{document}

\title{Minimizing atomic configurations \\ of short range pair potentials in two dimensions: \\
crystallization in the Wulff shape}

\author{Yuen Au Yeung$^1$, Gero Friesecke$^{1,2}$ and Bernd Schmidt$^1$ \\
{\small $^1$Zentrum Mathematik, Technische Universit\"at M\"unchen, Germany} \\[-1mm]
{\small $^2$Mathematics Institute, University of Warwick, U.K.}}

\date{September 4, 2009}

\maketitle

\begin{abstract}
We investigate ground state configurations of atomic systems in two dimensions
interacting via short range pair potentials.
As the number of particles tends to infinity, we show that low-energy configurations
converge to a macroscopic cluster of finite surface area and constant density, the latter
being given by the density of atoms per unit volume in the triangular lattice. 
In the special case of the Heitmann-Radin sticky disc potential and exact ground states, 
we show that the macroscopic cluster has a (unique) Wulff shape. This is done by showing that the 
atomistic energy, after subtracting off a bulk part and re-scaling, Gamma-converges
to a macroscopic anisotropic surface energy. 
\end{abstract}

\section{Introduction}
The question why a large number of atoms, under many conditions,
assembles into a periodic, crystalline structure of special
geometric shape remains poorly understood from a microscopic point of view.
Ultimately this should be a derivable consequence of the
underlying interatomic interactions, governed by the laws of quantum mechanics. 

In this paper, we address this question in the case of zero temperature,
simplified interatomic interactions (namely short range pair potentials), and
two instead of three dimensions. We believe that our methods will remain useful
in more general settings, including three-dimensional problems. In particular,
the prototype case allows to understand the key aspects of (1) formation of a local
lattice structure (2) emergence, as the number $N$ of particles gets large,
of a well defined anisotropic ``surface energy'' contribution of order
$O(N^{d/(d-1)})$ to the atomistic energy of an (exact or approximate)
atomistic minimizer (3) emergence of an overall geometric (Wulff) shape of
the atomistic minimizer as a consequence of surface energy minimization. 

We now state our results precisely. Consider $N$ particles in two dimensions, with
positions $x_1,...,x_N\in\R^2$, whose interaction energy is of pair potential
form,
\be \label{energy}
     E(x_1,\ldots,x_N) = \sum_{i\neq j} V(|x_i-x_j|).
\ee
We will assume that the atomic interaction potential is minimal when $|x_i-x_j|=1$
and short range; see conditions 
(H1), (H2) and (H3) in the next section for a precise statement. One then expects, at
least when the potential well of $V$ is sufficiently deep and narrow, the following
two phenomena to occur: 

(i) $E$ has crystallized ground states for any $N$.

(ii) In the limit $N\to\infty$, ground state configurations should assume, in a 
sense to be clarified, a particular overall geometric shape. 
\\[2mm]
By (i) we mean the following:
\\[1mm]
{\bf Definition} 
{\it We say that an energy $E \, : \, (\R^2)^N\to\R$ has crystallized ground states
if its infimum is attained and any minimizer -- after translation and rotation -- 
is a subset of the triangular lattice} 
\be \label{lattice} \vspace*{-2mm}
    \mcL := \{ m e_1 + n e_2 : m, \, n\in\Z\}, \;\;
    e_1=\begin{pmatrix} 1 \\ 0 \end{pmatrix}, \;\;
    e_2=\mbox{$\frac12$}\begin{pmatrix} 1 \\ \sqrt{3} \end{pmatrix}.
\ee
The energy $E$ given by (\ref{energy}) is rigorously known \cite{HeitmannRadin:80, Radin:81} 
to have crystallized ground states when $V$ is given by the 
Heitmann-Radin `sticky disc' potential
\be \label{radinpotential}
     V(r) = \left\{\begin{array}{ll}   +\infty, & 0\le r < 1 \\
                                       -1,      & r=1 \\
                                       0,       & r>1,\end{array}\right.
\ee
or the `soft disc' potential discussed in \cite{Radin:81}. For analogous results for the 
closely related sphere-packing problem in two dimensions see
\cite{Harborth, Conway, Toth}, and for insights into more general potentials see
\cite{Theil, E}. 

The present paper focuses on question (ii), which has remained open even in the case of energies
such as (\ref{radinpotential}). To describe the overall shape
of optimal subsets mathematically, we use the following strategy.

(1) associate to any atomic configuration
$\{x_1,\ldots,x_N\}$ its empirical measure $\sum_{i}\delta_{x_i}$, 

(2) re-scale it so that the total mass and the expected diameter of the support of 
a minimizing configuration remain of order one as
$N\to\infty$, 

(3) pass to the limit $N\to\infty$, 

(4) prove that the limit measure is a constant multiple of
a characteristic function of a set of finite perimeter (i.e., a characteristic
function belonging to the space $BV(\R^2)$), the constant being given by 
the density of atoms per unit volume in the triangular lattice

(5) derive, from atomistic energy minimization, a continuum Wulff-Herring type variational 
principle for the shape, via Gamma-convergence,

(6) use exact solubility of Wulff-Herring variational principles (cf. Taylor \cite{Taylor2}, 
Fonseca-M\"uller \cite{FonsecaMueller}) to identify the shape. 

The idea not to try and parametrize an atomistic configuration $\{x_1,..,x_N\}$
by displacements from a reference configuration (Lagrangian viewpoint), but to study its
empirical measure $\sum_{i=1}^N\delta_{x_i}$ (Eulerian viewpoint) and identify the limit measure
via a Gamma-convergence result, was recently introduced into the study of many-particle energies
by Capet and Friesecke \cite{CapetFriesecke}, in the context of Coulomb systems. 
Note that passage to the empirical measure eliminates the indefiniteness of $E$ under 
particle relabelling. 

The first part of our program, (1)--(4), has been carried out without needing to know whether $E$ has
crystallized ground states (in fact, in doing so one needs to prove a local form of crystallization).
The last two steps, (5)--(6), have been achieved {\it assuming} crystallization. 
This leads to the following two theorems.

The first result is not limited to minimizers, but applies to arbitrary states
whose energy difference from the ground state is of order $O(N^{1/2})$. 
\begin{theorem} \label{T1} (Formation of clusters with constant density and finite perimeter)
Suppose the energy $E$ is given by (\ref{energy}) and the interatomic
potential satisfies (H1), (H2), (H3). Let $\{x_1^{(N)},\ldots,x_N^{(N)}\}$ be
any sequence of connected (see Definition 3.1 below) $N$-particle configurations satisfying an energy bound of form
$$
     E(\{x_1^{(N)},\ldots,x_N^{(N)}\})\le -6N +  C N^{1/2}
$$
for some constant $C$ independent of $N$. Let $\{\mu_N\}$ be the associated sequence of
re-scaled empirical measures
\be \label{C:radonmeasure}
    \mu_N = \frac{1}{N} \sum_{i=1}^N \delta_{\frac{x_i^{(N)}}{\sqrt{N}}}.
\ee
Then: \\[1mm]
(i) Up to translation (that is to say, up to replacing $\mu_N$ by
$\mu_N(\cdot \,  + a_N)$ for some $a_N\in\R^2$) and passage to a subsequence, 
$\mu_N$ converges weak* in $\calM(\R^2)$ to $\mu\in\calM(\R^2)$. \\[2mm]
(ii) The limit measure is of the form
$$
   \mu = \rho\chi_{E},
$$
where $\rho=2/\sqrt{3}$ (i.e., the density of atoms per unit volume of the triangular lattice
$\mcL$) and $E$ is a set of finite perimeter of volume $1/\rho$.
\end{theorem}
In fact, any set $E$ of finite perimeter and volume $1/\rho$ can occur in the limit, as we prove in Section 5.

\begin{figure}[h!]
\begin{center}
\includegraphics[height=10cm]{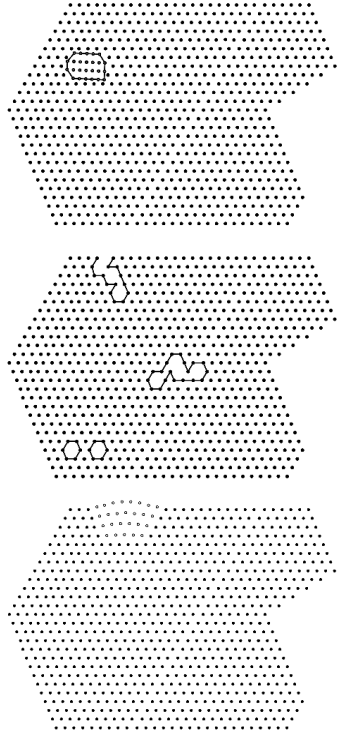}
\vspace*{-5mm}

\end{center}
\caption{\label{F:atomistic_configurations} Different atomistic configurations satisfying the assumptions of Thm.\ref{T1}.}
\end{figure}

Note also that on the atomistic level, quite irregular configurations are admitted by our
hypotheses. For instance the approximating atomistic configurations may contain elastic deformations,
cracks and vacancies, or inclusions of phases with different lattice structure, as long as these
only occupy regions of lengthscales smaller
than $N^{1/4}$, see Figure \ref{F:atomistic_configurations}.
Theorem \ref{T1} says that on the macro-scale, there nevertheless
result well defined clusters of constant, crystalline density. 

\noindent

For exact minimizers the ensuing cluster has a unique shape:

\begin{theorem} \label{Intro:T1} 
Suppose the energy $E$ is given by (\ref{energy}) and the interatomic
potential satisfies (H1), (H2), (H3). Assume in addition that $E$ has crystallized ground states (as is rigorously known
e.g. when $V$ is given by \eqref{radinpotential}). 
Let $\{x_1^{(N)},\ldots,x_N^{(N)}\}$ be any minimizing $N$-particle configuration of $E$, and let $\mu_N$ be the associated
re-scaled empirical measure (\ref{C:radonmeasure}).
As $N\to\infty$, up to translation and rotation (that is to say, up to replacing $\mu_N$ by
$\mu_N(R_N \, \cdot \,  + a_N)$ for some rotation $R_N\in SO(2)$ and some translation vector
$a_N\in\R^2$) $\mu_N$ converges weak star to the limit measure
\be \label{limitmeasure}
    \mu = \frac{2}{\sqrt{3}} \chi_{\overline{h}}
\ee
where $\overline{h}$ is the regular hexagon $\operatorname{conv}\, \{\pm e_1, \pm e_2, \pm (e_2-e_1)\}$.
\end{theorem}

\begin{figure}[h!]
\begin{center}
\includegraphics[height=4cm]{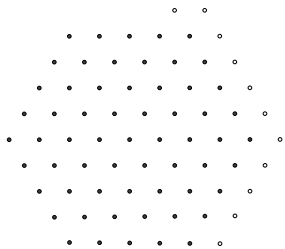}
\end{center}
\vspace*{-10mm}

\caption{\label{F:hex} A ground state of (\ref{energy}), (\ref{radinpotential}) for $N=72$.}
\end{figure}

Finally we remark that macroscopic uniqueness of the limit shape in Theorem \ref{Intro:T1} contrasts with an
unexpectedly large amount of non-uniqueness of the discrete minimizers. In a companion paper, we prove that the optimal
bound on the difference between two 
ground state configurations $\mu_N$ and $\mu_N'$,  
$$ \min \{\|\mu_N' - \mu_N(R\cdot + a)\| : R \in O(2), a \in \R^2 \} $$ 
(in suitable norms) scales like $N^{-1/4}$. Here optimal means that there exists a sequence $N_j\to\infty$ of particle
numbers for which there is a matching lower bound. 
Note that simple rearrangements of surface atoms 
only lead to differences of order $N^{-\frac{1}{2}}$ in the associated empirical measures.
Hence the $N^{-1/4}$ law shows unexpectedly strong fluctuations of finite ground state configurations
about the limiting Wulff shape.

\section{Atomistic energy}\label{sec:atomisticenergy}
Our object of study are
low-energy states of many-particle potential energy functionals of the form
\be \label{energy2}
     E(x_1,\ldots,x_N) = \sum_{i\neq j} V(|x_i-x_j|)
\ee
on $(\R^2)^N$, where $x_1,\ldots,x_N\in\R^2$ are the particle positions and the potential $V$
is assumed to satisfy the following hypotheses (see Figure \ref{F:soft}):

\begin{figure}[htbp]
\begin{center}
\includegraphics[height=3.7cm]{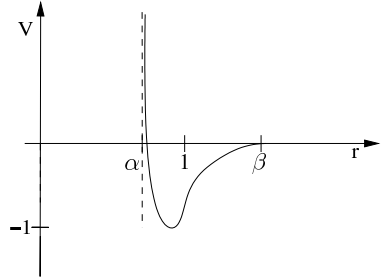}
\caption{\label{F:soft} Soft Heitmann-Radin potential.}
\end{center}
\end{figure}

\vspace{-5mm}
\begin{itemize}
\item[(H1)] (minimum at $r=1$) $V(1)=-1$, $V(r)>-1$ for all $r\neq 1$
\item[(H2)] (behavior at short and long range) There exist constants $\alpha\in(0,1]$,
$\beta\in [1,\infty)$ such that $V(r)=+\infty$ for $r<\alpha$, $V(r)=0$ for $r>\beta$, $V$
continuous on $(\alpha,\beta)$
\item[(H3)] (narrow potential well) The constants $\alpha$, $\beta$ from (H2) satisfy the
condition that the ball of radius $\beta$ contains at most six points whose distance from
the center and mutual distance is $\ge\alpha$.
\end{itemize}
Simple geometric considerations show that (H3) is always satisfied when $\alpha=1-\epsilon$,
$\beta=1+\epsilon$, and $\epsilon>0$ sufficiently small. The Heitmann-Radin sticky potential
(\ref{radinpotential}) is contained as a special case ($\alpha=\beta=1$).

The above hypotheses
are not aimed at maximum generality, but at simplicity of proofs. An important feature is
that unlike (\ref{radinpotential}), they allow elastic deformations.

As a simple consequence of the above hypotheses, one has the following lower bound
on the ground state energy:
\be \label{lowerbound}
    \inf_{x_1,\ldots,x_N\in\R^2} E(x_1,\ldots,x_N)  \ge - 6N.
\ee
This is because by (H2), (H3), each particle can have a negative interaction energy
with at most 6 other particles, and by (H1), the interaction energy with each of them is $\ge -1$.

Simple trial configurations show that the true ground state energy differs from the lower bound by
at most $O(N^{1/2})$. More precisely we claim that
\be \label{upperbound}
    \inf_{x_1,\ldots,x_N\in\R^2} E(x_1,\ldots,x_N) \le - 6N + a N^{1/2} + b
\ee
where we may take $a=4\sqrt{3}$, $b=12$. Indeed, consider trial states which are hexagonal subsets
of the triangular lattice:
$$
   \{x_1,\ldots,x_N\} = \Omega \cap \calL, \;\; \calL \mbox{ as in }(\ref{lattice}), \;\;
h_R\subset\Omega\subset\overline{h_R},
$$
where $\overline{h_R}$ is the closed hexagon of radius $R>0$ with ``bottom edge'' parallel to $e_1$,
\be \label{hexagon}
   \overline{h_R} = \operatorname{conv}\, \{ \pm Re_1, \, \pm R_{e_2}, \, \pm R(e_2\!-\! e_1)\},
\ee
$h_R$ is its interior, and $R$ and $\Omega$ are chosen suitably such that $\# \Omega\cap\calL = N$.
To infer (\ref{upperbound}),
it suffices to estimate the energy of these states from above by the right hand side of (\ref{upperbound}).
For completeness we include the (elementary) argument.

First, consider the case when $\Omega=\overline{h_R}$. In this case the number $N_s$ of ``surface atoms'',
i.e. the number of atoms in $\partial h_R$, equals the length of $\partial h_R$, i.e.
\be \label{surfaceatoms}
                      N_s=6R.
\ee
The total number of atoms is obtained by summation over the number of atoms in $\partial h_r$ for
$r\le R$,
\be \label{totalatoms}
   N = 1 + \sum_{r=1}^R 6r = 3R^2 + 3R + 1.
\ee
Solving equations (\ref{surfaceatoms}), (\ref{totalatoms}) for $N_s$ in terms of $N$ yields
$N_s = \sqrt{12 N - 3}-3$. To infer the energy, we only need to count the number of neighbors
of each atom, where $x$ is called a neighbor of $y$ if $|x-y|=1$.
Since the 6 ``corner'' atoms have three missing neighbors, and
the remaining $N_s-6$ surface atoms have two missing neighbors, the energy is
\be \label{etrial}
   E(x_1,\ldots,x_N) = -6N + 6 \cdot 3 + (N_s-6) \cdot 2 = -6N + 2 \sqrt{12 N - 3}.
\ee
Now add a partial layer of $k$ atoms in $\partial h_{R+1}$ (see Figure ?), where $k$ ranges from
$1$ to $6(R+1)-1$.
Denoting the positions of the layer atoms by $x_{N+1},\ldots,x_{N+k}$, the energy of the configuration is
\begin{eqnarray}
    E(x_1,\ldots,x_N,x_{N+1},\ldots,x_{N+k}) & = & E(x_1,\ldots,x_N) \\
                                     & + & \sum_{\ell=1}^k
    \Bigl(- 2 \cdot \# (\mbox{neighbors of $x_{N+\ell}$ in $\overline{h_R}$}) \\
                                     &  & \;\;\;\; - 1 \cdot \#
    (\mbox{neighbors of $x_\ell$ in $\partial h_{R+1}$})\Bigr). \label{ewithlayer}
\end{eqnarray}
Here the factor 2 appears because interactions between $\overline{h_R}$ and
the layer $\partial h_{R+1}$ appear only once in the sum over $\ell$, while
the layer-layer interactions appear twice.
If $k=1$, it is immediate from (\ref{etrial}) that $E(x_1,\ldots,x_N, x_{N+1})$ is bounded
from above by the right hand side of (\ref{upperbound}), so let us assume $k\ge 2$.
We may arrange the
layer so that there are at most 5 ``corner'' atoms with only 3 neighbors, and 2 ``end''
atoms with only 3 neighbors, while all other atoms in the layer have 4 neighbors (see Figure \ref{F:hex}).

The 5 corner atoms have 1 neighbor in $\overline{h_R}$
and 2 in $\partial h_{R+1}$, the end atoms have 2 neighbors in $\overline{h_R}$ and
1 in $\partial h_{R+1}$, and the remaining layer atoms have 2 neighbors in each set. Consequently
by (\ref{ewithlayer}) and (\ref{etrial})
\begin{eqnarray*}
    E(x_1,\ldots,x_N,x_{N+1},\ldots,x_{N+k}) & \le & E(x_1,\ldots,x_N)
                                  + 5(-2\cdot 1 - 1\cdot 2) \\
          & & + 2(-2\cdot 2 - 1\cdot 1)
    + (k-7)(-2\cdot 2 - 1\cdot 2) \\
                                     & = & -6(N+k) + 2\sqrt{12 N - 3} + 12.
\end{eqnarray*}
Estimating the square root trivially by $\sqrt{12(N+k)}$ yields the desired upper bound by the right
hand side of (\ref{upperbound}). Although this was not needed here, we remark that by the results of 
\cite{HeitmannRadin:80} this trial configuration 
actually forms a ground state of the potential \eqref{radinpotential}.

\section{Compactness and mass conservation} \label{S:Compactness}

The results in this section are not limited to minimizers, but apply to arbitrary states
in which the energy difference from the ground state is of order $O(N^{1/2})$. See Figure 
\ref{F:atomistic_configurations}.

Also, we confine ourselves to connected atomic configurations (see the Definition below).
In case of disconnected configurations, our analysis can be applied separately to the connected
components. Note also that minimizers must always be connected.
\\[2mm]
{\bf Definition 3.1} A finite set $S\subset\R^2$ of particle positions is
called connected if for any two $x, \, y\in S$ there exist $x_0,\ldots,x_N\in S$ such that $x_0=x$, $x_N=y$,
and the distance between successive points $x_{j-1}$, $x_j$ lies within the interaction range of the potential,
i.e. $|x_j-x_{j-1}| \le \beta$ for all $j=1,\ldots,N$.
\\[2mm]
Also, in the sequel we use the following standard notation. $C_0(\R^2)$ denotes the space of
continuous functions on $\R^2$ such that $f(x)\to 0$ as $|x|\to\infty$, $\calM(\R^2)$ denotes
the space of Radon measures on $\R^2$ of finite mass (recall $\calM$ is the dual of $C_0$),
and a sequence of Radon measures $\mu_N$ is said to converge weak* to $\mu$, notation:
$\mu_N\weakstar\mu$, if $\int_{\R^2}f\, d\mu_N\to \int_{\R^2}f\, d\mu$ for all $f\in C_0(\R^2)$.
\begin{proposition} \label{P1} (Compactness and mass conservation) Assume that the interatomic
potential satisfies (H1), (H2), (H3). Let $\{x_1^{(N)},\ldots,x_N^{(N)}\}$ be
any sequence of connected $N$-particle configurations satisfying the energy bound
$$
     E(\{x_1^{(N)},\ldots,x_N^{(N)}\})\le -6N +  C N^{1/2}
$$
for some constant $C$ independent of $N$. Let $\{\mu_N\}$ be the associated sequence of Radon
measures \eqref{C:radonmeasure}.
Then up to translation (that is to say, up to replacing $\mu_N$ by
$\mu_N(\cdot \,  + a_N)$ for some $a_N\in\R^2$) there exists
a subsequence converging weak* in $\calM(\R^2)$ to $\mu\in\calM(\R^2)$. Moreover the limit
measure satisfies $\mu_N\ge 0$, $\int_{\R^2}d\mu = 1$.
\end{proposition}
\Proof Throughout the proof, we write $\calS=\{x_1^{(N)},\ldots,x_N^{(N)}\}$, and denote
by $C$ a constant independent of $N$ whose value may change from line to line.

Since the $\mu_N$ are nonnegative and have mass $1$, they are bounded in
${\calM}(\R^2)$ and hence, by the Banach-Alaoglu theorem, there exists a weak* convergent
subsequence. Clearly the limit $\mu$ is nonnegative. It remains to show the only really
nontrivial assertion above, namely mass conservation $\int_{\R^2}d\mu=1$.

The key is to show that
\be \label{distbound}
    \mbox{diam}\, \calS = \max_{x,y\in \calS} |x-y| \le C N^{1/2}.
\ee
If such a bound holds, then after translation there exists a fixed ball of radius $R$ such that
$\mbox{supp}\,\mu_N\subset B_R$ for all $N$. Choose an increasing sequence of functions $\phi_n\in C_0(\R^2)$
such that $0\le\phi_n\le 1$, $\phi_n=1$ on $B_{nR}$. By dominated convergence of $\phi_n$ to $1$,
the weak* convergence of $\mu_N$ to $\mu$, and the fact that $\mbox{supp}\,\mu_N\subset B_R$,
$$
   \int d\mu = \lim_{n\to\infty}\int \phi_n\, d\mu = \lim_{n\to\infty}
   (\lim_{N\to\infty}\int \phi_n \, d\mu_N) = \lim_{n\to\infty} 1 = 1,
$$
completing the proof of the proposition.

It remains to establish the bound (\ref{distbound}). We begin by introducing a notion of
{\it neighbors} and a notion of {\it local energy}. Let $\epsilon\in(0,1)$, small
enough so that condition (H3) is satisfied for $\alpha = 1-\epsilon$ and $\beta = 1+\epsilon$. 
\\[2mm]
{\bf Definition} We say that $x\in S_N$ is a neighbor of $y\in S_N$ if $|x-y|\in[1-\epsilon, \,
1+\epsilon]$. The set of neighbors of $x\in \calS$ is denoted by $\calN(x)$.
\\[2mm]
Also, independently of the above notion of neighbors, we define a local energy,
as follows:
$$
    \Eloc(x) := \sum_{y\in \calS\backslash\{x\}} V(|x-y|),
$$
so that
$$
    E(\calS) = \sum_{x\in\calS} \Eloc(x).
$$
By hypotheses (H1)--(H3) on the potential $V$ and the finiteness of the energy of $\calS_N$,
it follows that
\be \label{Nrigid}
    \#\calN(x)\le 6, \;\; \Eloc(x)\ge -6,
\ee
with both bounds being sharp as is seen by considering points $x$ in the lattice $\calL$. A key point
now is that when the number of neighbors of $x$ is not equal to $6$, the local energy of $x$ is
bounded away from its optimum by a finite amount:
\be \label{ebound}
    \#\calN(x)<6 \Longrightarrow \Eloc(x)\ge -6 + \Delta,
    \mbox{ with }\Delta := \min_{r\not\in(1-\epsilon,1+\epsilon)} V(r) - \min_{r}V(r) >0.
\ee
Next, we construct an appropriate set in the plane associated with the configuration $\calS$.
For each $x\in\calS$, let $\calV(x)$ be the Voronoi cell of $x$,
\be \label{voronoi}
   \calV(x) = \Bigl\{y\in\R^2 : |y-x| \le |y-x'| \mbox{ for all }x'\in\calS\backslash\{x\}\Bigr\}.
\ee
As $\calV(x)$ may be unbounded, it is useful to introduce in addition its truncation
\be \label{truncvoronoi}
   \calV_{trunc}(x) := \calV(x) \cap B_1(x),
\ee
where $B_r(x)$ is the ball $\{y\in\R^2 : |y-x|\le 1\}$. We then define
\be \label{Omega}
   \Omega := \bigcup_{x\in\calS} \calV_{trunc}(x).
\ee
Elementary geometric considerations show that when $\epsilon>0$ is chosen sufficiently small, then
\be \label{noboundary}
    \#\calN(x)=6 \Longrightarrow
    \partial\Omega\cap \partial \calV_{trunc}(x) = \emptyset,
\ee
that is to say the cells associated to points with the maximum number of neighbors do not contribute
to the boundary of $\Omega$. Finally, since by construction $\calV_{trunc}(x)\subset B_1(x)$ and
$\calV_{trunc}(x)$ is convex, we have the following bound on the length of its boundary
\be \label{lengthbound}
    |\partial \calV_{trunc}(x)|\le 2\pi.
\ee
Consequently, denoting
\be \label{boundary}
     \partial \calS := \{x\in\calS : \#\calN(x)<6\},
\ee
using the plausible fact proved in Lemma \ref{L:con} below that
due to the connectedness of $\calS$ the set $\Omega$ is connected, and (\ref{noboundary}),
(\ref{lengthbound})
\be \label{bound1}
    \diam\calS \;\le\; \diam\Omega \;\le\; \mbox{$\frac12$} |\partial\Omega| \;\le\; \mbox{$\frac12$} \!\!\!\! \sum_{ x\in\calS\, | \,
    \calV_{trunc}(x)\cap\partial\Omega\neq\emptyset } \!\!\!\! |\partial\calV_{trunc}|
    \;\le\; \pi \,\#\partial\calS.
\ee
On the other hand, recalling $N=\#\calS$ and using the assumption on $E(\calS)$ in the lemma
and (\ref{ebound}),
\be
   -6N + CN^{1/2} \ge E(\calS) \ge -6N + \Delta\,\#\partial\calS,
\ee
and consequently
\be \label{bound2}
      \#\partial\calS \le \frac{C}{\Delta} N^{1/2}.
\ee
Combining (\ref{bound1}), (\ref{bound2}) establishes (\ref{distbound}).
The proof of the proposition is complete, except for the
following elementary geometric lemma. \eop
\begin{lemma} \label{L:con} If a configuration $\calS$ is connected in the sense of the
definition in Section \ref{S:Compactness} and the energy $E(\calS)$ is finite,
then the set $\Omega$ defined in (\ref{Omega}) is connected.
\end{lemma}
\Proof We begin by establishing an elementary inequality relating the constants $\alpha$ and $\beta$
appearing in the hypotheses on the potential $V$. By (H3), the ball of radius $\beta$ around the origin
does not contain the regular heptagon around the origin with sidelength $\alpha$, or equivalently
$$
    \frac{\alpha}{2} > \beta\sin\frac{\phi}{2}, \;\;\; \phi=\frac{2\pi}{7}.
$$
Hence
\be \label{betabound}
   \beta < \frac{\alpha}{2\sin\phi/2} = 1.152382... \, \alpha.
\ee
To prove the lemma, it suffices to show that if $|x-y|<\beta$, the truncated Voronoi cells around $x$ and $y$
have nonzero intersection. To establish this, it is enough to show that the midpoint $m=\frac{x+y}{2}$
belongs to $B_1(x)\cap B_1(y)$ and $V(x)\cap V(y)$. The first inclusion is immediate from
$|m-x|<\beta/2$, (\ref{betabound}) and $\alpha\le 1$. The second inclusion is equivalent to
$$
   |m-z| > \frac{|x-y|}{2} \; (=|m-x|=|m-y|)
$$
for all $z\in\calS\backslash\{x,y\}$. Now any such $z$ belongs to
$\R^2\backslash(B_\alpha(x)\cup B_\alpha(y)$ (otherwise the energy of the configuration would be infinite),
and any closest point $p$ to $m$ in the latter set has distance $|p-m|=\sqrt{\alpha^2 - (|x-y|/2)^2}$,
so it suffices to show that $\sqrt{\alpha^2 - (|x-y|/2)^2}>|x-y|/2$, or equivalently
\be \label{alphabound}
    \alpha > \frac{|x-y|}{\sqrt{2}}.
\ee
But thanks to (\ref{betabound}) we have
$$
   |x-y|<\beta< 1.152382... \, \alpha < \sqrt{2}\alpha.
$$
This establishes (\ref{alphabound}), completing the proof of the lemma. \eop

\section{Proof of formation of clusters with constant density and finite perimeter}\label{sec:clusters}

We now prove Theorem \ref{T1}, i.e. we show that under the assumptions of Proposition \ref{P1}, the limit measure $\mu$ of
the re-scaled empirical measures of atomistic configurations is a constant multiple of
a characteristic function., the constant being given by the density of atoms per unit volume of the triangular lattice.
\\[2mm]
{\bf Proof of Theorem \ref{T1}. } As in the proof of Proposition \ref{P1} we use the decomposition of $\R^2$ into truncated Voronoi
cells $V_{trunc}(x)$ (see (\ref{voronoi})--(\ref{truncvoronoi})) associated with the points $x$ of an atomistic configuration 
${\cal S}$. The main technical idea of the
proof is to introduce and investigate the following {\it volume excess function}
\be \label{volexcess}
  \Phi(x) := \Bigl|  \frac{1}{|V_{trunc}(x)|} - \frac{1}{|h'_{1/\sqrt{3}}(x)|} \Bigr|.
\ee
Here $h'_r(x)$ is the closed hexagon of radius $r$ around $x$ which is obtained from $\overline{h_r}$ (see \eqref{hexagon}) by
a 30$^o$ rotation around the centre. This hexagon is the Voronoi cell of any interior point $x$ of a particle configuration
${\cal S}$ on the triangular lattice \eqref{lattice}. 

Besides the sequence of re-scaled empirical measures (\ref{C:radonmeasure}), we will make use of the following
auxiliary sequences:
\begin{eqnarray}
  \mutildeN & = & \frac1N \sum_{x\in\calS_N}  \frac{\chi_{N^{-1/2}V_{trunc}(x)}}{|N^{-1/2}V_{trunc}(x)|}, \label{eq:mutN} \\
  \muttildeN &= & \frac1N \sum_{x\in\calS_N}
\frac{\chi_{N^{-1/2}V_{trunc}(x)}}{|N^{-1/2}h'_{1/\sqrt{3}}(x)|}. \label{eq:muttN}
\end{eqnarray}
{\bf Step 1} First we claim that
\be \label{firstdiff}
    \mu_N - \mutildeN \weakstar 0 \mbox{ in }\calM(\R^2).
\ee
Indeed, since for any $x_0\in\calS_N$, $V_{trunc}(x_0)\subseteq B_1(x_0)$, we have for all test functions
$\phi\in C_0(\R^2)$
\begin{eqnarray*}
  \lefteqn{\Bigl| \int_{\R^2} \left( \delta_{x_0/\sqrt{N}} 
  - \frac{\chi_{N^{-1/2}V_{trunc}(x_0)}}{|N^{-1/2}V_{trunc}(x_0)|} \right)
   \phi \Bigr| } \, \\
  & = & \Bigl| \frac{1}{|N^{-1/2}V_{trunc}(x_0)|} \int_{N^{-1/2}V_{trunc}(x_0)} \Bigl(\phi(x_0/\sqrt{N}) - \phi(x) \Bigr)
       \, dx \Bigr| \\
  & \le & \sup_{|x-y|\le N^{-1/2}} |\phi(x)-\phi(y)|
\end{eqnarray*}
and consequently
\begin{eqnarray*}
   \Bigl| \int \phi \, d\mu_N - \int \phi \, d\mutildeN\Bigr| 
   & \le & \frac{1}{N}\sum_{x\in\calS_N} \sup_{|x-y|\le N^{-1/2}} |\phi(x)-\phi(y)| \\ 
   &  =  & \sup_{|x-y|\le N^{-1/2}} |\phi(x)-\phi(y)|  \to 0 \quad (N\to\infty),
\end{eqnarray*}
establishing (\ref{firstdiff}).\smallskip 

\noindent 
{\bf Step 2} Next, we compare the measures $\mutildeN$ and $\muttildeN$. We will show that
\be \label{seconddiff}
    \muttildeN - \mutildeN \to 0 \mbox{ in }L^1(\R^2)
\ee
(and hence, a forteriori, $\muttildeN - \mutildeN \weakstar 0$ in $\calM(\R^2)$).
First, we decompose both measures into their ``interior'' and ``boundary'' parts, as follows. Here
$\partial\calS_N$ is as defined in (\ref{boundary}), and $\operatorname{int} \calS_N := \calS_N\backslash\partial\calS_N$.
$$
   \mutildeN = \lambdatildeN + \nutildeN, \; \lambdatildeN = \frac1N \sum_{x\in \operatorname{int}\calS_N}
                                               \frac{\chi_{N^{-1/2}V_{trunc}(x)}}{|N^{-1/2}V_{trunc}(x)|}, \;
   \nutildeN = \frac1N \sum_{x\in \partial S_N}
                                               \frac{\chi_{N^{-1/2}V_{trunc}(x)}}{|N^{-1/2}V_{trunc}(x)|}
$$
and analogously for $\muttildeN$. Roughly speaking, we will argue that the difference between $\nuttildeN$ and
$\nutildeN$ is small because they are small separately, due to the fact that the number of boundary atoms grows only
like $N^{1/2}$, and that the difference between $\lambdattildeN$ and $\lambdatildeN$ is small because otherwise
this would cost elastic energy.

To make the first argument precise, we use that for any configuration
$\calS=\{x_1,\ldots,x_N\}$ with finite energy, $V_{trunc}(x_i)\supseteq B_{\alpha/2}(x_i)$, due to the fact
that $|x_j-x_i| \ge \alpha$ for all $j\neq i$, so that the nearest point in $\calS$ to $y\in B_{\alpha/2}(x_i)$ is $x_i$. Hence
$$
   0 \le |\nuttildeN - \nutildeN| \le \sum_{x\in\partial\calS} \chi_{N^{-1/2}V_{trunc}(x)} \Bigl(
     \underbrace{\frac{1}{|B_{\alpha/2}|} + \frac{1}{|h'_{1/\sqrt{3}}|}}_{=:K}\Bigr).
$$
Since $\#\partial\calS_N\le (C/\Delta)N^{1/2}$ (see (\ref{bound2})) and
$V_{trunc}(x)\subseteq B_1(x)$, we infer that
\be \label{nubound}
   \| \nuttildeN - \nutildeN \|_{L^1(\R^2)} \le \frac{C\, K}{\Delta} N^{1/2}
\max_{x\in\calS_N}\|\chi_{N^{-1/2}V_{trunc}(x)}\|_{L^1}
   \le \frac{C\, K}{\Delta} N^{-1/2} |B_1| \, \to \, 0 \; (N\to\infty).
\ee

Next we analyze the difference between the interior parts,
$$
   \lambdatildeN - \lambdattildeN = \sum_{x\in \operatorname{int}\calS_N} \chi_{N^{-1/2}V_{trunc}(x)} \Bigl(
   \frac{1}{|V_{trunc}(x)|} - \frac{1}{|h_{1/\sqrt{3}}(x)|}\Bigr).
$$
For any point $x\in \operatorname{int}\calS_N$, $\#\calN(x)=6$, so $\calN(x)=\{y_1,\ldots,y_6\}$,
where we may assume that the $y_j$ are numbered so that $y_j-x=r_j(\cos\phi_j,\sin\phi_j)$, $0\le\phi_1<\phi_2<...<\phi_6<2\pi$.
For elementary geometric reasons, namely that the only way to arrange the $y_i$ so that
$|y_j-x|=1$ for all $j$ and $|y_j-y_{j-1}|=1$ for all $j$ is to place them at the corners of a regular hexagon around $x$,
we have
$$
  |V_{trunc}(x)|\to |h_{1/\sqrt{3}}(x)| \; \mbox{ if }
  \sum_{y\in N(x)} V(|x-y|) + \sum_{j=1}^6 V(|y_{j+1}-y_j|) \to -12.
$$
Hence by continuity, given $\delta>0$ there exists $\Delta(\delta)>0$ such that for all $x\in \operatorname{int}\calS_N$ the
volume excess (\ref{volexcess}) satisfies the following implication: 
\begin{eqnarray}
   & &  \Phi(x) \ge \delta \nonumber \\
   & &  \Longrightarrow  \tilde{E}_{\ell oc}(x) :=
        \sum_{y\in \calN(x)} V(|x-y|) + \sum_{j=1}^6 V(|y_{j+1}-y_j|) \ge -12 + \Delta. \label{nextbound}
\end{eqnarray}
It will be convenient to extend the function $\tilde{E}_{\ell oc}(x)$ to all of $\calS_N$, by setting it equal to $-12$
when $x\in\partial\calS_N$.
Combining the energy bound assumed in Theorem \ref{T1}, the fact that each particle pair appears in $\tilde{E}_{\ell oc}(x)$
for at most four $x\in\calS_N$ whereas it appears twice in $E$, and (\ref{nextbound}) yields
\begin{eqnarray*}
   - 6N + C\, N^{1/2} & \ge & E(\calS) \; \ge \; \frac12\sum_{x\in\calS} \tilde{E}_{\ell oc}(x) \\
   & \ge & \frac12 \Bigl( - 12 \# \Bigl\{ x\in\calS_N : x\in\partial\calS_N \mbox{ or }\Phi(x)<\delta\Bigr\} \\ 
   & &  + (-12+\Delta) \# \Bigl\{x\in \operatorname{int}\calS_N : \Phi(x) \ge \delta\Bigr\}\Bigr) \\
   & \ge & -6N + \frac{\Delta}{2} \# \Bigl\{x\in \operatorname{int}\calS_N : \Phi(x) \ge \delta\Bigr\}.
\end{eqnarray*}
Consequently
$$
   \#\Bigl\{ x\in \operatorname{int} \calS_N : \Phi(x)\ge\delta\Bigr\} \le \frac{2\, C}{\Delta} N^{1/2}
$$
and so
\begin{eqnarray*}
  \|\lambdattildeN - \lambdatildeN \|_{L^1(\R^2)} & = &
  \int_{\R^2} \Bigl| \sum_{x\in \operatorname{int} \calS_N}
  \chi_{N^{-1/2}V_{trunc}(x)} \Bigl(
   \frac{1}{|V_{trunc}(x)|} - \frac{1}{|h_{1/\sqrt{3}}(x)|}\Bigr) \Bigr| \\
  & \le & \sum_{x\in \operatorname{int} \calS_N} \Phi(x) \|\chi_{N^{-1/2}V_{trunc}(x)}\|_{L^1} \\
  & \le & N^{-1} |B_1| \sum_{x\in \operatorname{int}\calS_N} \Phi(x) \\
  & \le & N^{-1} |B_1| \Bigl( \delta \, N + \underbrace{\max_{x\in \operatorname{int}\calS_N}\Phi(x)}_{\le K} \;
     \underbrace{\#\Bigl\{x\in \operatorname{int}\calS : \Phi(x)\ge\delta\Bigr\} }_{\le \frac{2C}{\Delta}N^{1/2}} \Bigr).
\end{eqnarray*}
Letting $N\to\infty$ gives
$$
   \limsup_{N\to\infty} \| \lambdattildeN - \lambdatildeN \|_{L^1(\R^2)} \le |B_1| \delta.
$$
Since $\delta>0$ was arbitrary, the left hand side equals zero, that is to say $\lambdattildeN - \lambdatildeN \to 0$
in $L^1(\R^2)$. Together with (\ref{nubound}), this establishes (\ref{seconddiff}).\smallskip 

\noindent 
{\bf Step 3} Having established that the limits of $\mu_N$, $\mutildeN$, $\muttildeN$ all coincide, it suffices to
study the limiting behavior of the sequence $\muttildeN$, which -- as we shall see -- is compact in a much ``stronger''
space. After a change on a set of measure zero,
$$
   \muttildeN = \frac{1}{|h'_{1/\sqrt{3}}|} \chi_{N^{-1/2}\Omega_N}, \;\;\; \Omega_N = \bigcup_{x\in\calS_N}V_{trunc}(x).
$$
Clearly the $\muttildeN$ only take values in $\{0,\rho\}$, belong to the space $BV(\R^2)$, are bounded in $L^1(\R^2)$, and -- by (\ref{distbound}) -- after
translation are supported in some fixed ball of radius $R$. We claim that they are also bounded in $BV(\R^2)$. This is because
$|\partial\Omega_N|$ is bounded by a constant times $N^{1/2}$, by (\ref{bound1}) and (\ref{bound2}). Consequently by the Banach-Alaoglu theorem
and the compact embedding $BV(B_R)\hookrightarrow L^1(B_R)$, a subsequence converges weak* in $BV$ and strongly in $L^1$ to some limit
$\tilde{\tilde{\mu}}\in BV$. By the strong $L^1$ convergence, $\tilde{\tilde{\mu}}$ only takes values in $\{0,\rho\}$, and by
Steps 1 and 2, $\tilde{\tilde{\mu}}=\mu$. The proof of Theorem \ref{T1} is complete. \eop

\section{Gamma-convergence of surface energy and emergence of Wulff shape}

In case of low-energy states, we now study the shape of the limiting
cluster $E$ obtained in the previous section.

\begin{figure}[htbp]
\begin{center}
\includegraphics[height=5cm]{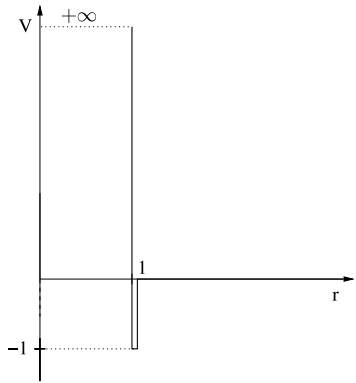}
\caption{\label{fig:stickyHR} The Heitmann-Radin `sticky disc' potential.}
\end{center}
\end{figure}

Our current methods are restricted to exact minimizers and to energies with crystallized ground states. 
Here the problem simplifies because by (H1), (H2) and (H3) we may without loss of generality assume 
that the interaction potential is given by the Heitmann-Radin potential (\ref{radinpotential}) (see Figure \ref{fig:stickyHR}).

For such configurations, we can derive a limiting variational principle, as follows. 
It is useful (as in \cite{CapetFriesecke}) to re-formulate the minimization problem for the atomistic energy $E$ in terms of empirical
measures instead of particle configurations. Define the following energy functionals on the set $\calA$
of probability measures on $\R^2$ (i.e., nonnegative Radon measures on $\R^2$ of mass 1):
\be \label{newenergy}
    I_N(\mu) := \begin{cases}
    \int_{\R^4\bs diag}N V(N^{1/2}|x-y|)\,d\mu\otimes d\mu, & \mu = \frac{1}{N}\sum_{i=1}^N \delta_{x_i/\sqrt{N}} \\ 
       & \mbox{ for some distinct }x_i\in \calL\\
    +\infty, & \mbox{otherwise}.
    \end{cases}
\ee
This definition says that $I_N(\mu_N)=E(x_1,\ldots,x_N)$ when $\mu_N$ is the re-scaled empirical measure 
(\ref{C:radonmeasure}) of the configuration $\{x_1,\ldots,x_N\}$. In particular, the re-scaled empirical measure
minimizes $I_N$ if and only if the underlying configuration minimizes $E$.

We now show:
\begin{theorem}\label{T2} The sequence of functionals $N^{-1/2}(I_N + 6N)$ Gamma-converges, with respect
to weak$^\ast$ convergence of probability measures, to the limit functional $I_\infty \, : \, \calA \to \R\cup\{\infty\}$
given by 
\be \label{limitfunctional}
   I_\infty(\mu) := \begin{cases}
   \int_{\p^\ast E}\Gamma(\nu_E) d\calH^1(x), 
   & \mu = \frac{2}{\sqrt{3}} \chi_E~ \mbox{\rm for some set $E$ of } \\ 
      & \mbox{ {\rm finite perimeter and mass} } \frac{\sqrt{3}}{2}, \\
   +\infty,& \mbox{{\rm otherwise}}, 
   \end{cases}
\ee
where $\Gamma$ is the function
\be \label{hexenergy}
   \Gamma(\nu) = 2 \left(\nu_2 - \frac{\nu_1}{\sqrt{3}} \right) \mbox{{\rm for }}\nu={-\sin\varphi \choose \cos\varphi}, \;
   \varphi\in[0,\mbox{$\frac{2\pi}{6}$}],
\ee
extended $\frac{2\pi}{6}$--periodically.
\end{theorem} 

\begin{figure}[htbp]
\begin{center}
\includegraphics[height=4cm]{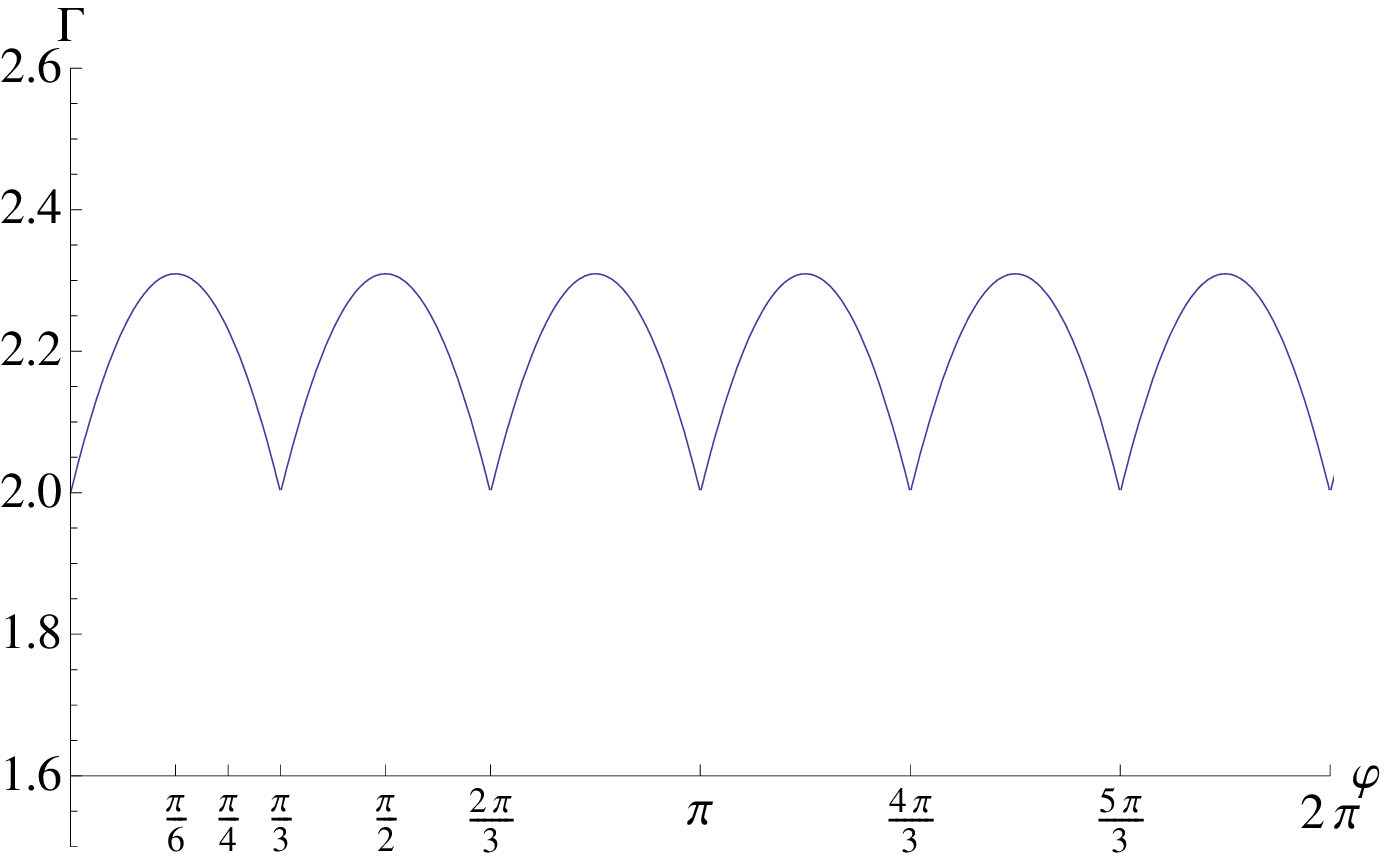}
\quad\includegraphics[height=4cm]{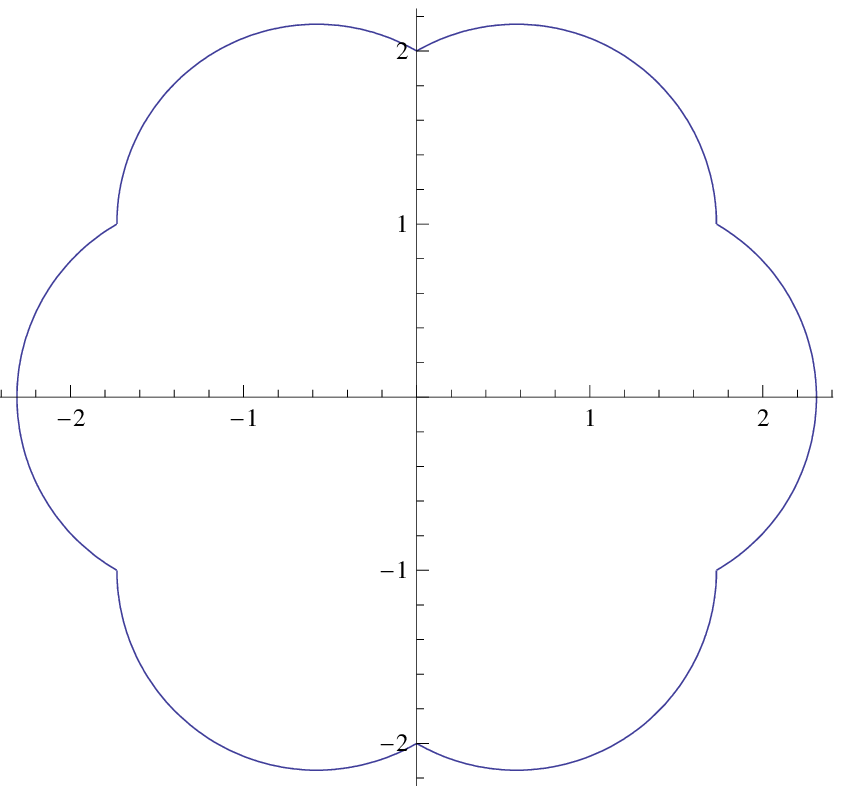}
\caption{Plot of the limiting surface energy density $\Gamma$ in dependence of $\varphi$ and polar plot of $\Gamma$ as a function of the normal $\nu$.}
\end{center}
\end{figure}

The physical idea which is made quantitative here that crystalline anisotropy gives rise to nontrivial surface energies goes 
back to Wulff (1901) and Herring (1951) in the physics literature. Rigorous interfacial energy results have a large
and sophisticated statistical mechanics literature, particularly in the context of the Ising model (see e.g. the book
by Dobrushin, Kotecky, and Shlosman \cite{DKS}). For a previous treatment of an interfacial energy problem via 
Gamma-convergence see \cite{ABC}. 

\Proof 
{\em Lower bound.} Let $\calS_N = \{x_1^{(N)},\ldots,x_N^{(N)}\}
\subset \calL$ be any sequence of $N$-particle configurations
whose associated sequence of Radon measures $\mu_N$ (see \eqref{C:radonmeasure}) 
weak$^\ast$-converges to a probability measure $\mu \in \calA$. We need to show that $\liminf I_N(\mu_N) \ge I_\infty(\mu)$.  

Associate to $\calS_N$ the following auxiliary set
$$ 
   H_N := \bigcup_{x \in \calS_N}
   \overline{N^{-1/2} h'_{1/\sqrt{3}}(x)},
$$
where $h'_r(x)$ denotes the open hexagon around $x$ of radius
$r$ introduced in the previous section. Note that  
for $x \in \calL$, the closure of $h'_{1/\sqrt{3}}(x)$ is the Voronoi cell of $x$ with 
respect to the complete lattice $\calL$. 

The boundary $\p H_N$ is a disjoint union of simple closed polygons $V_1, 
\ldots, V_M$. 

\begin{figure}[htbp]
\begin{center}
\includegraphics[height=6.2cm]{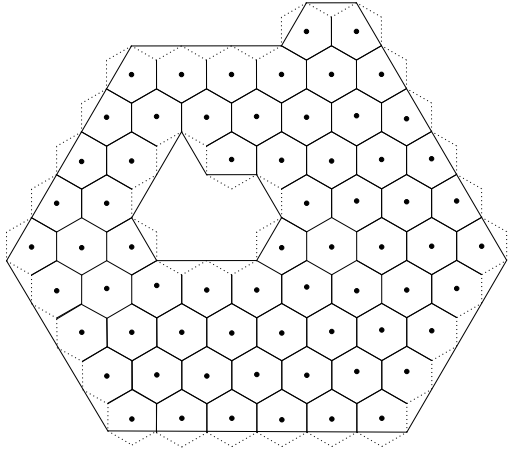}
\caption{Constructing $H_N'$ from $H_N$.}
\label{F:HNprime}
\end{center}
\end{figure}

Because the boundary $\p H_N$ oscillates on the atomic scale, we 
define yet another auxiliary set $H_N'$ (see Figure \ref{F:HNprime}) which removes these oscillations and will hence
allow us to obtain a sharp lower bound on $I_N(\mu_N)$ via standard weak lower semicontinuity results
on surface functionals. 
If $V_j \subset \p H_N$ 
is a simple closed polygon, then $V_j = \bigcup_{i = 1}^m 
[v_{i+1}, v_i]$ with segments $[v_{i + 1}, v_i]$ of length $1 / \sqrt{3N}$ for 
$$ v_1, \ldots, v_m \in \left(\frac{1}{3\sqrt{N}}(e_1 + e_2) + \frac{1}{\sqrt{N}}\calL\right) \cup 
   \left(\frac{1}{3\sqrt{N}}(2 e_2 - e_1) + \frac{1}{\sqrt{N}}\calL\right), \quad v_{m + 1} := v_1. $$
As the corner points $v_i$ alternate between the lattices $\frac{1}{3\sqrt{N}} 
(e_1 + e_2) + \frac{1}{\sqrt{N}}\calL$ and $\frac{1}{3\sqrt{N}}(2 e_2 - e_1) 
+ \frac{1}{\sqrt{N}}\calL$ comprising the 
dual lattice of $\calL$, $m$ is an even number, and $V_j' = 
\bigcup_{i = 1}^{m/2} [v_{2i-1}, v_{2i+1}]$ is a simple closed polygon. We 
now define the set $H_N' \subset \R^2$ as the unique closed set with $\calS_N 
\subset H_N'$ such that $\p H_N' = \bigcup_{J = 1}^M V_j'$. It can be easily 
seen that
\be\label{eq:HNHNstrich}
  |H_N' \triangle H_N| 
   \le \frac{\#\p \calS_N}{8N\sqrt{3}}.
\ee
Here, for two sets $A$ and $B$, $A\triangle B$ denotes the symmetric difference $(A\backslash B)\cup (B\backslash A)$. 
Note that $[v,w] \cap [x,y] \ne \emptyset$ defines a one-to-one correspondence 
between boundary segments $[v,w]$ of $H_N$ and nearest neighbor bonds $[x,y]$ of 
$\frac{1}{\sqrt{N}}\calL$, i.e.\ $x, y \in \frac{1}{\sqrt{N}}\calL$ with $|x - y| 
= \frac{1}{\sqrt{N}}$, such that $x \in \frac{1}{\sqrt{N}}\calS_N$ and $y \notin 
\frac{1}{\sqrt{N}}\calS_N$. But then 
there is a one-to-two correspondence between boundary segments $[v,w]$ of $H_N'$ 
and nearest neighbour bonds $[x,y]$ such that $x \in \calS$ and $y \notin \calS$. 
Since boundary segments of $H_N'$ are of length $\frac{1}{\sqrt{N}}$, this shows 
that 
\be\label{eq:HNstrichEnergie}
  I_N(\mu_N) + 6N = 2 \sqrt{N} \calH^1(\p H_N'). 
\ee

Precisely as in Step 1 in the proof of Theorem \ref{T1} we see that 
\be\label{eq:muprimeN}
  \mu_N' := \frac{1}{N} \sum_{x \in \calS_N}
  \frac{\chi_{N^{-1/2}h'_{1/\sqrt{3}}}(x)}{|N^{-1/2}h'_{1/\sqrt{3}}(x)|}
\ee
weak$^\ast$ converges to $\mu$. 
Now---similarly as in Step 3 in  the proof of 
Theorem \ref{T1}---we obtain that after a change on a set of measure zero 
$$ 
   \mu_N' = \rho \chi_{H_N}. 
$$
Moreover, if $\frac{1}{\sqrt{N}} (I_N(\mu_N) + 6N)$ is bounded it follows from \eqref{eq:HNstrichEnergie}
that the restrictions $\mu_N|_{B_R}$ are bounded in $BV(B_R)$ for each $R>0$, and hence converge
to $\mu|_{B_R}$ in $L^1(B_R)$, due to the compact embedding $BV(B_R)\hookrightarrow L^1(B_R)$.
Since $\|\mu\|_{L^1(\R^2)} = \|\mu_N\|_{L^1(\R^2)} = 1$ for all $N$, we obtain 
that even $\mu_N \to \mu$ in $L^1(\R^2)$ and $\mu = \rho \chi_E$ for a set $E$ 
of finite perimeter and area $\frac{1}{\rho}$. 
 
Furthermore, if (a subsequence of) $\frac{1}{\sqrt{N}} (I_N(\mu_N) + 6N)$ 
is bounded, then by \eqref{eq:HNHNstrich} and \eqref{bound2} 
$\rho \chi_{H_N'} \to \mu = \rho \chi_E$ in $L^1(\R^2)$, whereas by \eqref{eq:HNstrichEnergie}, 
\eqref{hexenergy} and the fact that each boundary segment of $H_N'$ is parallel to 
$e_1$ , $e_2$ or $e_1-e_2$ 
$$ \frac{1}{\sqrt{N}} (I_N(\mu_N) + 6N) 
   = 2 \calH^1(\p H_N') 
   = \int_{\p(H_N')} \Gamma(\nu)\,d\mathcal{H}^1.
$$
Now since $\Gamma$, extended to a $1$-homogeneous function on $\R^n$, is convex and 
satisfies a growth condition of the type $\Gamma(\nu) \geq c |\nu|$ for some $c > 0$, 
a lower semicontinuity result for SBV functions (cf.\ Theorem 5.22 in Ambrosio-Fusco-Pallara
\cite{Ambrosio-Fusco-Pallara}) establishes the lower bound 
$$ \liminf_{N \to \infty} \frac{1}{\sqrt{N}} (I_N(\mu_N) + 6N) 
   \ge \int_{\p^\ast E}\Gamma(\nu_E) d\calH^1(x). 
$$
\smallskip

\noindent
{\em Upper bound.} The idea to find a recovery sequence is to
approximate gradually a set of finite perimeter $E$ by sets of
simpler geometric shape. To be more precise, a set of
finite perimeter will be approximated by $C^\infty$ sets, $C^\infty$
sets by sets with polygonal boundaries and polygonal sets by
polygonal sets 
having all their corners in $\frac{1}{n}\calL$ for some $n \in \N$. 
The recovery sequence for $E$ will then be extracted through a
diagonalization process. Finally, a suitable continuity property of surface
integrals,
$$ 
   \lim_{N \to \infty} \int_{\p^\ast P_N} \Gamma(\nu(x))\,d\calH^1(x)
   = \int_{\p^\ast P} \Gamma(\nu(x))\,d\calH^1(x),
$$
where $(P_N)$ is the approximation sequence and $P$ denotes the
approximated set in the respective step, will complete the proof. 
\smallskip

\noindent 
{\bf Step 1} First, let $P \subset \R^2$ be a bounded set with polygonal boundary such 
that every corner of $\p P$ lies in $\frac{1}{n} \calL$ for some $n \in \N$. Assume that 
the volume of $P$ is $\frac{\sqrt{3}}{2} + \alpha_n$. Consider the sequence 
of particle configurations $\tilde{\calS}_{n,N} = \calL \cap \sqrt{N} P$ consisting 
of $M = M_{n,N}$ atoms. Elementary geometric considerations show that 
\be\label{eq:diffMN} 
   |M_{n, N} - N| \le c ( \alpha_n N + \sqrt{N} ).  
\ee 
for a constant $c$ independent of $n$ and $N$. Let $\mu_{n,N}$ denote the associated 
rescaled empirical measures. Clearly, $\mu_{n,N} \weakstar \rho \chi_P$ as $N \to \infty$. 
In addition, a straightforward calculation of the associated energy by evaluating the 
surface energy contribution 
$$ 
   \frac{1}{\sqrt{N}} \sum_{x \in \p \tilde{\calS}_{n,N}} 
   \#\left\{ x' \in \calL \setminus \tilde{\calS}_{n,N} : |x-x'| = 1
   \right\},
$$ 
along line segments of $\p P$ gives  
\be \label{step1-upperbound}
  \left| \frac{1}{\sqrt{N}} ( I_M(\mu_{n,N}) + 6M ) 
  -  \int_{\p P} \Gamma(\nu(x)) d\,\calH^1(x) \right| 
  \le c \left( \alpha_n + \frac{1}{\sqrt{N}} \right)
\ee
for a constant $c$ independent of $n$ and $N$. 
Indeed, an elementary argument shows that if $S$ is a boundary segment of $P$ of 
length $L$ with normal $\nu$, then the number of pairs $(x,y) \in \tilde{\calS}_{n,N} 
\times (\calL \setminus \tilde{\calS}_{n,N})$ with $|x - y| = 1$ such that the 
segment $[x, y]$ intersects $\sqrt{N} S$ is equal to 
$$ 
   \lfloor \sqrt{N} \rfloor \Gamma(\nu) L + O\left(L\right) 
   =  \sqrt{N} \left( \Gamma(\nu) L + O\left( \frac{L}{\sqrt{N}} \right) \right).  
$$
\smallskip 

\noindent 
{\bf Step 2} Now let $P \subset \R^2$ be any set of volume $\frac{\sqrt{3}}{2}$ with polygonal 
boundary. By perturbing the corners of $\p P$ slightly, it is easy to see that there is a 
sequence of polygonal sets $P_n$ whose boundary vertices lie in $\frac{1}{n}\calL$ satisfying 
\be\label{eq:PnDeltaP} 
   |P_n \triangle P| \to 0 \mbox{ as } n \to \infty  
\ee
and 
\be\label{eq:RandPnP} 
   \int_{\p^\ast P_n} \Gamma(\nu_{P_n}) \to \int_{\p^\ast P} \Gamma(\nu_{P}) \mbox{ as } n \to \infty.
\ee
Now choosing $n = n(N) \to \infty$ appropriately, we obtain a sequence of configurations 
$$ 
   \tilde{\calS}_{N} := \tilde{\calS}_{n(N), N} 
$$ from the configurations constructed in Step 1 such that the associated rescaled empirical measures 
$\tilde{\mu}_N$ satisfy $\tilde{\mu}_N \weakstar \rho \chi_P$ and 
$$ 
   \frac{1}{\sqrt{N}} ( I_M(\tilde{\mu}_N) + 6M ) 
   \to \int_{\p P} \Gamma(\nu(x)) d\,\calH^1(x)  
$$
by \eqref{eq:PnDeltaP}, \eqref{step1-upperbound} and \eqref{eq:RandPnP}.  

Now if $\# \tilde{\calS}_N \neq N$, then the configuration $\tilde{\calS}_N$ needs to
be modified either by adding or removing some set of points with negligible
surface term. By \eqref{eq:diffMN} $\frac{1}{N}|M - N| \to 0$ as $N \to \infty$. 
If $\tilde{\calS}_N$ contains less than $N$ particles, then we can add the missing number of points 
on lattice sites within some parallelogram whose sidelengths are of order $\sqrt{N - M}$ and 
such that the points in the parallelogram are a distance greater than $1$ away from $\tilde{\calS}_N$. 
Then none of the new points interacts with $\tilde{\calS}_N$. 
In case $M > N$, just remove $M-N$ points of $\tilde{\calS}_N$ lying in a common parallelogram 
whose sidelengths are of order $\sqrt{M - N}$. (This is always possible if $N$ is large enough.) 

Denoting the corresponding rescaled empirical measure by $\mu_N$ and recalling that 
$\frac{1}{N}|M - N| \to 0$ it is not hard to see that $\mu_N \weakstar \rho \chi_P$ as $N \to \infty$. 
Since furthermore 
$$ 
   \left| \frac{1}{\sqrt{N}} ( I_M(\tilde{\mu}_{N}) + 6M ) 
   -  \frac{1}{\sqrt{N}} ( I_N(\mu_{N}) + 6N ) \right|
   \le c \sqrt{\frac{1}{N}(M - N)} 
$$
for some constant $c > 0$, we also obtain 
$$ 
   \frac{1}{\sqrt{N}} ( I_N(\mu_N) + 6N ) 
   \to \int_{\p P} \Gamma(\nu(x)) d\,\calH^1(x)  
$$
by \eqref{eq:diffMN} and \eqref{step1-upperbound}. 
\smallskip

\noindent
{\bf Step 3} 
A recovery sequence for a general set $E$ of finite perimeter is now obtained by a diagonalization 
argument due to the following density and continuity results. 

Suppose first $E$ is a bounded set of finite perimeter with $C^{\infty}$-boundary and volume 
$\frac{\sqrt{3}}{2}$. By piecewise linear approximations of $\p E$ and scaling 
with factors close to $1$ we easily construct approximations $P_n$ which have 
the same volume, a polygonal boundary and satisfy 
$$ 
   \chi_{P_n} \weakstar \chi_E 
   \quad\mbox{ and }\quad 
   \int_{\p^\ast P_n} \Gamma(\nu_{P_n}) \to \int_{\p^\ast E} \Gamma(\nu_{E}). 
$$

Now let $E$ be a bounded set of finite perimeter with volume $\frac{\sqrt{3}}{2}$. As shown e.g.\ in 
\cite[Proposition 4.7 and Remark 4.8]{Braides:98}, there are bounded sets $E_n$ of finite perimeter with $C^{\infty}$-boundary and---
after rescaling with a factor close to $1$---volume $\frac{\sqrt{3}}{2}$ such that 
$$ 
   \chi_{E_n} \weakstar \chi_E \mbox{ and } 
   \int_{\p^\ast E_n} \Gamma(\nu_{E_n}) \to \int_{\p^\ast E} \Gamma(\nu_{E}). 
$$

Finally note that a truncation argument yields that an analogous result holds for approximating sets 
of finite perimeter with bounded sets of finite perimeter. This concludes the proof. \eop


To complete the proof of Therem \ref{Intro:T1}, it remains to infer convergence of minimizers. A technical
detail we need to pay attention to is that, unlike in many other Gamma-convergence results, here sequences with bounded
energy are not in general compact in the topology in which the Gamma-convergence occurs. This is because, due to the
translation invariance of the functionals $N^{-1/2}(I_N+6N)$, sequences can lose mass as $N\to\infty$. 
\\[2mm]
{\bf Proof of Theorem \ref{Intro:T1}} Let $\{x_1^{(N)},\ldots,x_N^{(N)}\}$ be any minimizing $N$-particle configuration of $E$, 
and let $\mu_N$ be the associated re-scaled empirical measure (\ref{C:radonmeasure}). By the connectedness of minimizing
configurations and Proposition \ref{P1}, after suitable translations $\mu_N \mapsto \mu_N(\cdot + a_N)$ the limit measure
$\mu$ has full mass. Hence sequences of exact minimizers are compact in the topology in which the
Gamma-convergence occurs (namely weak* convergence of probability measures). By standard arguments in Gamma-convergence,
$\mu$ is a minimizer of the limit functional $I_\infty$. 

\begin{figure}[htbp]
\begin{center}
\includegraphics[height=5cm]{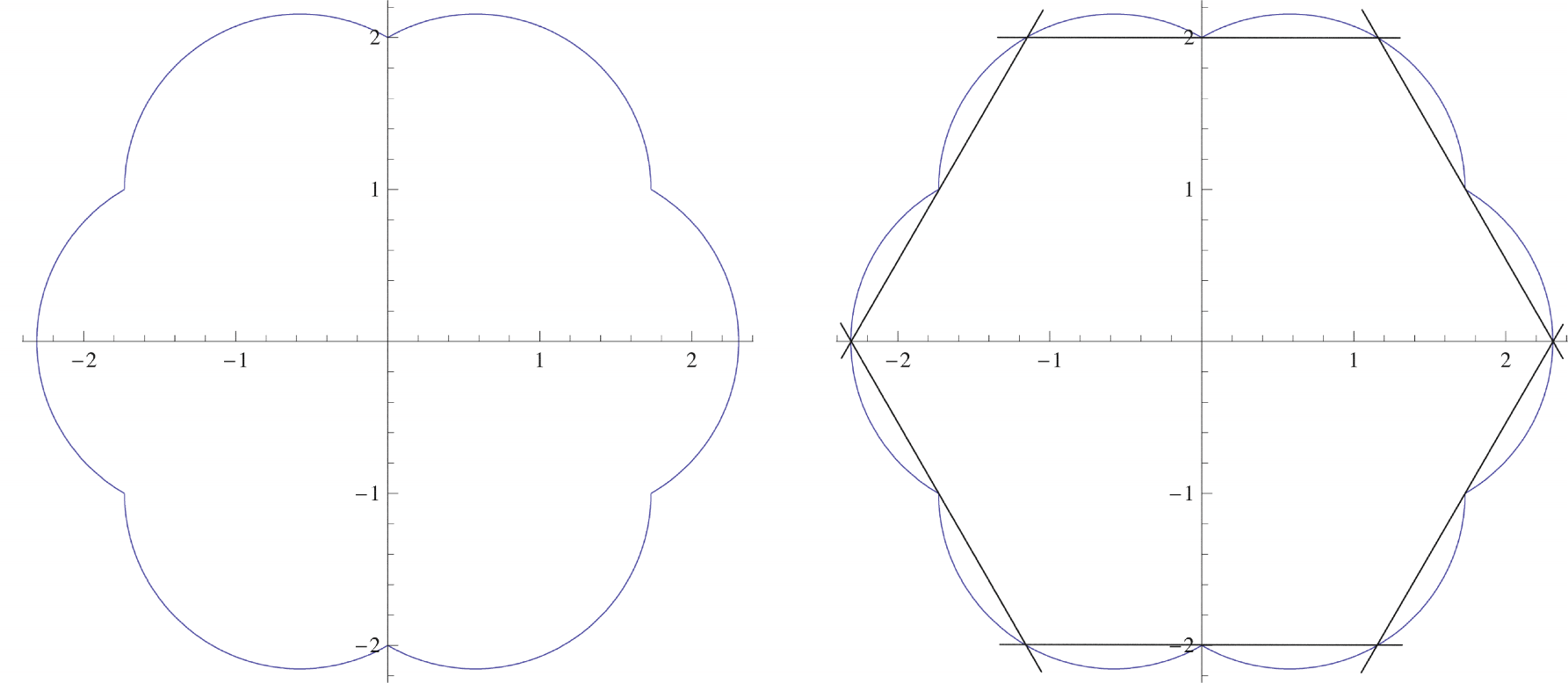}
\caption{\label{fig:wulff} Polar plot of $\nu \to \Gamma(\nu)$ and the associated Wulff shape $W_{\Gamma}$.}
\end{center}
\end{figure}

We now appeal to the uniqueness theorem for Herring type
energies due to Taylor (in the language of geometric measure theory), in the version by
Fonseca and M\"uller (who work in the present setting of boundary integrals for sets of finite perimeter):
\begin{theorem}\label{theo:TaylorFonsecaMueller} {\rm \cite{Taylor2, FonsecaMueller}} A functional of form
$$
     I(E) = \int_{\partial^\ast E} \Gamma(\nu(x))\, d{\cal H}^{n-1}(x),
$$
with $\Gamma \, : \, S^{n-1}\to[0,\infty)$ continuous and bounded away from zero,
is minimized over sets $E\subset\R^n$ of finite perimeter and volume $1$ if and only if $E$ agrees,
up to translation and up to a set of measure zero, with $\lambda W_\Gamma$, where $W_\Gamma$ is
the Wulff set
$$
      W_\Gamma := \Bigl\{x\in\R^n : x\cdot \nu\le \Gamma(\nu) \mbox{ for all }\nu\in S^{n-1}\Bigr\}
$$
and $\lambda>0$ is the unique normalization constant such that $\lambda W_\Gamma$ has volume 1.
\end{theorem}
In the present case of the energy (\ref{hexenergy}), an elementary calculation shows that the
Wulff set is given by the intersection of the six half-spaces $x\cdot\nu \le \Gamma(\nu)$ for
the minimizing normals $\nu_{2\pi j/6}$, $j=1,\ldots,6$, i.e. a regular hexagon with bottom face
parallel to the $e_1$-axis (see Figure \ref{fig:wulff}). This completes the proof of Theorem \ref{Intro:T1}.

\bibliographystyle{alpha}

\vspace*{5mm}

\noindent
\centerline{\small Email of authors: {\tt yuen@online.de, gf@ma.tum.de, schmidt@ma.tum.de}}
\end{document}